# ASYMPTOTIC THEORY FOR THE SEMIPARAMETRIC ACCELERATED FAILURE TIME MODEL WITH MISSING DATA

By Bin Nan,[1] John D. Kalbfleisch and Menggang Yu

*University of Michigan, University of Michigan and Indiana University*

We consider a class of doubly weighted rank-based estimating methods for the transformation (or accelerated failure time) model with missing data as arise, for example, in case-cohort studies. The weights considered may not be predictable as required in a martingale stochastic process formulation. We treat the general problem as a semiparametric estimating equation problem and provide proofs of asymptotic properties for the weighted estimators, with either true weights or estimated weights, by using empirical process theory where martingale theory may fail. Simulations show that the outcome-dependent weighted method works well for finite samples in case-cohort studies and improves efficiency compared to methods based on predictable weights. Further, it is seen that the method is even more efficient when estimated weights are used, as is commonly the case in the missing data literature. The Gehan censored data Wilcoxon weights are found to be surprisingly efficient in a wide class of problems.

**1. Introduction.** Instead of modeling the hazard function for censored survival data, as in the Cox model [6], modeling the (transformed) failure time directly is sometimes appealing to practitioners since it postulates a simple relationship between the response variable and covariates with easily interpretable parameters. Let $T$ denote the failure time transformed by a known monotone function $h$, $C$ be the corresponding transformed censoring time, $\Delta = 1(T \leq C)$ and $Y = \min(T, C)$, where $1(\cdot)$ denotes an indicator function. The model of interest is

$$(1.1) \qquad T_i = \theta_0' Z_i + e_i, \qquad i = 1, \ldots, n,$$

Received January 2008; revised July 2008.
[1]Supported in part by NSF Grant DMS-07-06700.

where the $e_i$'s are independent and identically distributed (i.i.d.) with unknown distribution $F$, and $e_i$ is independent of $(Z_i, C_i)$ for all $i$. When $h = \log$, the model is called the accelerated failure time model (see, e.g., [12]).

For a cohort of $n$ i.i.d. observations of $X_i = (Y_i, \Delta_i, Z_i)$, $i = 1, \ldots, n$, [4] proposed an imputation type of least squares method, where the censored survival time is replaced by an estimate of the mean residual life conditional on the covariates, which is obtained from the Kaplan–Meier estimator on the residual scale. Stute [24, 25] proposed a weighted least squares method with weights obtained from the Kaplan–Meier estimator for the transformed survival time. [21, 26] and [30], among others, studied the rank-based estimating method and proved the asymptotic properties using martingale theory for counting processes.

In this article, we consider a general rank-based estimating method for model (1.1) in the presence of missing data as arise, for example, in case-cohort studies (e.g., [19, 23]) where data are missing by design. Specifically, let $Z_i = (Z'_{1i}, Z'_{2i})'$ and assume that $Z_{1i}$ is missing at random (see [14]), while $Z_{2i}$, $Y_i$ and $\Delta_i$ are always observed for all $i$. The situations where $Z_i = Z_{1i}$ for all $i$, or where $Z_{2i}$ is not included in model (1.1), are special cases. In the latter of these special cases, $Z_{2i}$ is usually called an auxiliary variable in the missing data literature. The approach in this article extends the work of [16] for case-cohort studies, where weights are predictable and the counting process approach of [26] applies. It can be applied to general two-phase outcome-dependent sampling designs for censored survival data and allows the use of nonpredictable weights that can yield more efficient parameter estimates. The proof of efficiency gains from using estimated weights, even though the true weights are given, similarly follows the approach of [18].

This article is organized as follows. In Section 2, we introduce the doubly weighted rank-based estimating method with arbitrary weights (i.e., either predictable or nonpredictable), and link the proposed estimating function to a semiparametric framework that is more suitable for applying empirical process theory. Methods based on both known weights and estimated weights are considered. We describe asymptotic properties of the proposed estimators in Section 3, with detailed proofs given in Section 6. In Section 4, we discuss the asymptotic efficiency and some simulation results that compare methods of using predictable weights and nonpredictable weights and methods of using known weights and estimated weights. We make a few concluding remarks in Section 5.

**2. Doubly weighted semiparametric estimating function.** For the $i$th subject, $Z_{2i}, Y_i$ and $\Delta_i$ are always observed. Let $R_i$ be the missing data indicator that takes value 1 if $Z_{1i}$ is also observed and 0 otherwise. Suppose that $Z_{1i}$ is missing at random, so that

$$\pi_i = \Pr(R_i = 1|Z_i, Y_i, \Delta_i) = \Pr(R_i = 1|Z_{2i}, Y_i, \Delta_i)$$

for each $i$. This holds, for example, when independent Bernoulli sampling is implemented in a two-phase sampling design that includes the case-cohort study as a special case.

To estimate $\theta_0$ in model (1.1), we follow [15] and define the following random map

$$\begin{aligned}\Psi_n(\theta, \eta, \rho) &= \frac{1}{n}\sum_{i=1}^{n}\psi(X_i; \theta, \eta, \rho) \\ &= \frac{1}{n}\sum_{i=1}^{n}\Omega_i \rho(Y_i - \theta'Z_i, \theta)\{Z_i - \eta(Y_i - \theta'Z_i, \theta)\}\Delta_i,\end{aligned}$$
(2.1)

where $\theta \in \Theta \subset \mathbb{R}^d$ is the $d$-dimensional Euclidean parameter of interest with unknown true value $\theta_0$, and $\eta$ and $\rho$ are real valued (vectors of) functions that can be viewed as infinite dimensional nuisance parameters.

When $\eta(t, \theta)$ is replaced by an estimator of the true function (see [21])

$$\eta_0(t, \theta) = \frac{E\{1(Y - \theta'Z \geq t)Z\}}{E\{1(Y - \theta'Z \geq t)\}},$$

with $\eta_0(t, \theta_0) = E(Z|Y - \theta_0'Z \geq t)$, random map (2.1) becomes a weighted estimating function for $\theta$, where $\Omega_i$ are subject specific weights and $\rho(t, \theta)$ is a weight function. Clearly such an estimating function is semiparametric.

To be more general, we assume that the true functional forms of $\eta$ and $\rho$ are unknown and need to be estimated, and study the estimating function $\Psi_n(\theta, \hat{\eta}_n, \hat{\rho}_n)$ with

$$\hat{\eta}_n(t, \theta) = \frac{\sum_{j=1}^{n} W_j 1(Y_j - \theta'Z_j \geq t)Z_j}{\sum_{j=1}^{n} W_j 1(Y_j - \theta'Z_j \geq t)},$$
(2.2)

where $W_j$ are subject specific weights that may or may not equal $\Omega_j$. This is the source of the term "double weights" (see [31]); the purpose of introducing two possibly different subject specific weights will soon become clear. A particularly interesting weight function $\rho(t, \theta)$ is taken to be $\rho_0(t, \theta) = \Pr(Y - \theta'Z \geq t)$, and it can be estimated by

$$\hat{\rho}_n(t, \theta) = \frac{\sum_{j=1}^{n} W_j 1(Y_j - \theta'Z_j \geq t)}{\sum_{j=1}^{n} W_j},$$
(2.3)

a weighted Gehan-type weight. This type of weight provides a very desirable property. The corresponding estimating function $\Psi_n(\theta, \hat{\eta}_n, \hat{\rho}_n)$ is monotone in $\theta$. See [31] for the detailed derivation.

In this article, we focus on the estimator of $\theta$ obtained from the estimating function $\Psi_n(\theta, \hat{\eta}_n, \hat{\rho}_n)$, where $\hat{\eta}_n$ is given in (2.2). The estimator $\hat{\rho}_n$ can be more flexible, but we will be particularly interested in the one given by (2.3). Using two possibly different sets of subject specific weights $\Omega_i$ and $W_i$ in $\Psi_n(\theta, \hat{\eta}_n, \hat{\rho}_n)$ yields great flexibility that covers a broad range of problems. The following are a few examples:

(i) When $\rho = 1$ and $\Omega_i = W_i = 1$ for all $i$, (2.2) becomes
$$\hat{\eta}_n(t, \theta) = \frac{\sum_{j=1}^n 1(Y_j - \theta' Z_j \geq t) Z_j}{\sum_{j=1}^n 1(Y_j - \theta' Z_j \geq t)},$$
and the estimating function $\Psi_n(\theta, \hat{\eta}_n, 1)$ becomes the rank-based estimating function studied by [26] and [29], among others. [26] and [30] proved asymptotic linearity of $\Psi_n(\theta, \hat{\eta}_n, 1)$ and thus normality of the estimator obtained from $\Psi_n(\theta, \hat{\eta}_n, 1) = 0$ using a stochastic integral formulation and martingale theory for counting processes.

(ii) When $\hat{\rho}_n$ takes the form in (2.3) and $\Omega_i = W_i = 1$ for all $i$, $\Psi_n(\theta, \hat{\eta}_n, \hat{\rho}_n)$ becomes the estimating function of [26] with Gehan weights. The monotonicity of such an estimating function was studied by [7].

(iii) When $\hat{\rho}_n$ takes the form in (2.3) and $\Omega_i = 1$, $W_i = 1(i \in \mathcal{SC})/\Pr(i \in \mathcal{SC})$ for all $i$ where $\mathcal{SC}$ denotes the set of labels of the subcohort in a case-cohort study, $\Psi_n(\theta, \hat{\eta}_n, \hat{\rho}_n)$ becomes the estimating function of [16] with generalized Gehan-type weights.

(iv) When $\hat{\rho}_n$ takes the form in (2.3) and $\Omega_i = 1$, $W_i = R_i/\pi_i$ for all $i$, where $\pi_i$ depends on $\Delta_i$, $\Psi_n(\theta, \hat{\eta}_n, \hat{\rho}_n)$ becomes an extension of the estimating function of [31] (where the authors focused on numerical aspects and did not provide asymptotic properties). The weights $\Omega_i = 1$ and $W_i = R_i/\pi_i$ have been applied to case-cohort studies to potentially improve efficiency in the Estimator II of [2] as well as in [5, 13] for the Cox model.

(v) When $\Omega_i = W_i = R_i/\pi_i$, the estimating function $\Psi_n(\theta, \hat{\eta}_n, \hat{\rho}_n)$ can be applied to a general missing data problem with covariate $Z_{1i}$ missing at random. This arises, for example, in a two-phase sampling design and yields an estimator that is similar to that proposed in [20] and further studied by [3] for the Cox model.

In examples (i), (ii) and (iii), the estimating functions can be formulated as martingales, and the related theory applies. In the last two situations, however, weights $\Omega_i$ and/or $W_i$ depend on $\Delta_i$, particularly in case-cohort studies, and, thus, are not predictable. There is no martingale representation of these weighted estimating functions. Further complications are: (1) the estimating function $\Psi_n(\theta, \hat{\eta}_n, \hat{\rho}_n)$ is a nonsmooth function of $\theta$, so that the methods developed for smooth estimating functions based on Taylor expansions do not apply; and (2) the nuisance parameters $\eta$ and $\rho$ are explicit functions of $\theta$, whereas usual semiparametric models assume that nuisance parameters do not vary with the parameter of interest.

Our simulation study shows a substantial efficiency gain when such outcome-dependent weights are used and more efficiency gain when the known weights are estimated from observed data. This latter result has been often noted (see, e.g., [3, 11, 18, 22], among many others). For these reasons, it is desirable to rigorously investigate the theoretical properties of the estimators obtained from the estimating function $\Psi_n(\theta, \hat{\eta}_n, \hat{\rho}_n)$ with both known and estimated flexible weights.

**3. Asymptotic properties.** Assume that the observed data are i.i.d. In addition to Conditions 1–3 in [30] (also assumed in [26]), we assume Conditions (A) and (B) below and derive asymptotic properties of the estimator obtained from the weighted estimating function $\Psi_n(\theta, \hat{\eta}_n, \hat{\rho}_n)$. In particular, these results apply when $\hat{\eta}_n$ is given by (2.2) and $\hat{\rho}_n$ takes the form (2.3), which estimates $\rho_0 = \Pr(Y - \theta'Z \geq t)$ with either true weights $W_i$ or their estimates $\hat{W}_i$. Our method does not depend on stochastic integrals and, hence, does not require predictability of the weights. So, it applies to a much broader range of estimating functions. Note that $\hat{\eta}_n(t, \theta)$ in (2.2) and $\hat{\rho}_n(t, \theta)$ in (2.3) are not differentiable in $\theta$.

CONDITION (A). There exist constants $\tau < \infty$ and $\xi$, such that $\Pr(Y - \theta'Z \geq \tau) \geq \xi > 0$ for all $Z$ and $\theta \in \Theta$.

CONDITION (B). The selection probability $\pi = \Pr(R = 1|Z_2, Y, \Delta) \geq \zeta > 0$ for all $Z_2$, $Y$ and $\Delta$ for some constant $\zeta$.

Condition (A) follows an assumption in equation (3.1) of [26]. Condition (B) is a common assumption in the missing data literature and guarantees that the inverse selection probability weights are bounded. Using empirical process theory, we follow the idea of [26] and [30] to show the asymptotic linearity of $\Psi_n(\theta, \hat{\eta}_n, \hat{\rho}_n)$ in $\theta$ in a neighborhood of the true value $\theta_0$. We adopt the empirical process notation of [27]. In particular, for a function $f$ of a random variable $U$ that follows distribution $P$, we define

$$Pf = \int f(u)\, dP(u),$$

$$\mathbb{P}_n f = n^{-1} \sum_{i=1}^{n} f(U_i),$$

$$\mathbb{G}_n f = n^{-1/2}(\mathbb{P}_n - P)f$$

and refer all the details to the reference. Throughout the article, we assume that $\Omega_i$ and $W_i$ are bounded and satisfy $E(\Omega_i|X_i) = E(W_i|X_i) = 1$, for all $i$, and set $\varepsilon_\theta = Y - \theta'Z$ and $\varepsilon_0 = Y - \theta_0'Z$.

3.1. *Using true weights.* Consistency and rate of convergence of the proposed estimator $\hat{\theta}_n$ for general $\eta$ and $\rho$ are given in Theorems 3.1 and 3.2, respectively. Asymptotic normality of $\hat{\theta}_n$ obtained from the estimating function $\Psi_n(\theta, \hat{\eta}_n, \hat{\rho}_n)$, with $\hat{\eta}_n$ and $\hat{\rho}_n$ taking the forms in (2.2) and (2.3), is given in Theorem 3.3. Proofs are deferred to Section 6.

THEOREM 3.1. *Denote $\Psi(\theta, \eta, \rho) = P[\rho(\varepsilon_\theta, \theta)\{Z - \eta(\varepsilon_\theta, \theta)\}\Delta]$. Let $\Theta$, the parameter space of $\theta$, be compact, assume that $\theta_0 \in \Theta$ is the unique solution of $\Psi(\theta, \eta_0, \rho_0) = 0$ and let $\|\cdot\|$ be the supremum norm. If $\|\eta - \eta_0\| \leq \delta_n$ and $\|\rho - \rho_0\| \leq \delta_n$ with $\delta_n \downarrow 0$, where $\eta$, $\eta_0$, $\rho$ and $\rho_0$ belong to Glivenko–Cantelli classes and are bounded, then:*

(i) *In outer probability,*

(3.1)  $$\|\Psi_n(\theta, \eta, \rho) - \Psi(\theta, \eta_0, \rho_0)\| \to 0;$$

(ii) *An approximate root $\hat{\theta}_n$ satisfying $\Psi_n(\hat{\theta}_n, \eta(\cdot, \hat{\theta}_n), \rho(\cdot, \hat{\theta}_n)) = o_{p^*}(1)$ is consistent;*

(iii) *When $\hat{\eta}_n$ and $\hat{\rho}_n$ are given respectively by (2.2) and (2.3), an approximate root $\hat{\theta}_n$ satisfying $\Psi_n(\hat{\theta}_n, \hat{\eta}_n(\cdot, \hat{\theta}_n), \hat{\rho}_n(\cdot, \hat{\theta}_n)) = o_{p^*}(1)$ is consistent.*

THEOREM 3.2. *Let $\Theta_0 \subset \Theta$ be a neighborhood of $\theta_0$, $\|\cdot\|_0$ be the supremum norm in $\Theta_0$ and $\hat{\eta}_n$ be as in (2.2). Assume that $\|\hat{\rho}_n - \rho_0\|_0 = O_{p^*}(n^{-1/2})$, and assume that both $\hat{\rho}_n$ and $\rho_0$ are bounded and belong to a Donsker class. Let $\hat{\theta}_n$ be an approximate root satisfying $\Psi_n(\hat{\theta}_n, \hat{\eta}_n(\cdot, \hat{\theta}_n), \hat{\rho}_n(\cdot, \hat{\theta}_n)) = o_{p^*}(n^{-1/2})$. Suppose $\Psi(\theta, \eta_0(\cdot, \theta), \rho_0(\cdot, \theta))$ is differentiable with bounded continuous derivative $\dot{\Psi}_\theta(\theta, \eta_0(\cdot, \theta), \rho_0(\cdot, \theta))$ in $\Theta_0$, and $\dot{\Psi}_\theta(\theta_0, \eta_0(\cdot, \theta_0), \rho_0(\cdot, \theta_0))$ is nonsingular. Then, $\|\hat{\eta}_n - \eta_0\|_0 = O_{p^*}(n^{-1/2})$ and $|\hat{\theta}_n - \theta_0| = O_{p^*}(n^{-1/2})$. Finally, if $\hat{\rho}_n$ takes the form in (2.3) and $\rho_0(t, \theta) = \Pr(\varepsilon_\theta \geq t)$, then the above conditions for $\hat{\rho}_n$ and $\rho_0$ are satisfied.*

In the proofs of the above theorems, given in Section 6, we apply the permanence of the Donsker property under closures and convex hulls (see [27]) to show that (2.2) and (2.3) and their limits are Donsker. A variety of sufficient conditions for Donsker classes of functions are provided in [27].

When $\hat{\eta}_n$ takes the form in (2.2), the estimating function $\Psi_n(\theta, \hat{\eta}_n, \hat{\rho}_n)$ is discontinuous in $\theta$. In the case of full cohort data with $\Omega_i = W_i = 1$ for all $i$, [21, 26, 30] showed, with considerable effort, the asymptotic linearity of $\Psi_n(\theta, \hat{\eta}_n, 1)$, in a neighborhood of the true parameter $\theta_0$, in order to prove asymptotic normality. [16] had equally complicated arguments for asymptotic linearity in case-cohort studies where the weights $W_i$ do not depend on $\Delta_i$. We avoid the stochastic integral formulation and apply empirical process theory to show the asymptotic linearity of $\Psi_n(\theta, \hat{\eta}_n(\cdot, \theta), \hat{\rho}_n(\cdot, \theta))$ around $\theta_0$ for the class of missing data problems considered here. In Theorem 3.3, we focus on the situation where $\hat{\eta}_n$ and $\hat{\rho}_n$ are, respectively, given by (2.2) and (2.3). For other types of bounded weight functions $\hat{\rho}_n$ and $\rho_0$, proofs of asymptotic normality follow the same steps, and the same asymptotic representation should hold if $\{\hat{\rho}_n\}$ and $\{\rho_0\}$ are Donsker and $\hat{\rho}_n$ is an asymptotic linear estimator. This approach takes care of both predictable and nonpredictable weights.

THEOREM 3.3. *Let $\hat{\eta}_n$ and $\hat{\rho}_n$ be as in (2.2) and (2.3). Let $\hat{\theta}_n$ be an approximate root satisfying $\Psi_n(\hat{\theta}_n, \hat{\eta}_n(\cdot, \hat{\theta}_n), \hat{\rho}_n(\cdot, \hat{\theta}_n)) = o_{p^*}(n^{-1/2})$. Let $\mathcal{Y}$ and $\mathcal{Z}$ denote the sample spaces of random variables $Y$ and $Z$, respectively. Suppose that $\rho_0(\varepsilon_\theta, \theta)$ and $\eta_0(\varepsilon_\theta, \theta)$ are differentiable in $\theta$ with derivatives $\dot{\rho}_{0\theta}$ and $\dot{\eta}_{0\theta}$, which are uniformly bounded and continuous in $\Theta_0 \times \mathcal{Y} \times \mathcal{Z}$. Note that this implies that $\Psi(\theta, \eta_0(\cdot, \theta), \rho_0(\cdot, \theta))$ is differentiable in $\theta$ with bounded continuous derivative $\dot{\Psi}_\theta(\theta, \eta_0(\cdot, \theta), \rho_0(\cdot, \theta))$ in $\Theta_0$. Then, we have the following:*

(i) *The asymptotic linearity*

$$
\begin{aligned}
(3.2) \quad & n^{1/2} \Psi_n(\hat{\theta}_n, \hat{\eta}_n(\cdot, \hat{\theta}_n), \hat{\rho}_n(\cdot, \hat{\theta}_n)) \\
&= n^{1/2} \Psi_n(\theta_0, \hat{\eta}_n(\cdot, \theta_0), \hat{\rho}_n(\cdot, \theta_0)) \\
&\quad + n^{1/2}(\hat{\theta}_n - \theta_0) \dot{\Psi}_\theta(\theta_0, \eta_0(\cdot, \theta_0), \rho_0(\cdot, \theta_0)) + o_{p^*}(1)
\end{aligned}
$$

*holds;*

(ii) *If $\dot{\Psi}_\theta(\theta_0, \eta_0(\cdot, \theta_0), \rho_0(\cdot, \theta_0))$ is nonsingular, then $n^{1/2}(\hat{\theta}_n - \theta)$ is asymptotically normal with the asymptotic representation*

$$
\begin{aligned}
(3.3) \quad n^{1/2}(\hat{\theta}_n - \theta_0) &= \{-\dot{\Psi}_\theta(\theta_0, \eta_0(\cdot, \theta_0), \rho_0(\cdot, \theta_0))\}^{-1} \\
&\quad \cdot \mathbb{G}_n \bigg[\Omega \rho_0(\varepsilon_0, \theta_0)\{Z - \eta_0(\varepsilon_0, \theta_0)\}\Delta \\
&\qquad - \int W \rho_0(t, \theta_0)\{Z - \eta_0(t, \theta_0)\} 1(\varepsilon_0 \geq t) \, d\Lambda_0(t)\bigg] \\
&\quad + o_{p^*}(1).
\end{aligned}
$$

REMARK. As becomes clear in the proof of Theorem 3.3, the asymptotic representation (3.3) is the same if the weight function $\rho_0(t, \theta)$ is known, and, in fact, such a property does not depend on what $\rho_0(t, \theta)$ is. This finding is consistent with the claim in Section 4 of [26]. Equation (3.3) reduces to the result of [16] for predictable $W$ when $\Omega = 1$ and $\rho_0(t, \theta) = 1$. The variance estimator for $\hat{\theta}_n$ can be obtained following the method described in [16] based on the asymptotic representation (3.3) and the original idea of [9]. Alternative variance estimation methods can be found in [10, 17]. Later, in Section 4.1, we show that letting $\Omega = W$ yields more efficient estimation for the example of a case-cohort study.

3.2. *Using estimated weights.* In Theorems 3.1, 3.2 and 3.3, the subject-specific weights $W_i$ and $\Omega_i$ are assumed to be known. This is a reasonable assumption for many types of sampling designs when weights are the inverse of sampling probabilities, because sampling probabilities are usually prespecified by investigators. In the missing data literature, many authors (e.g., [22] and [3]) have pointed out that using the estimated weights improves the asymptotic efficiency, even though the true weights are known. Suppose true weights $W_i$ are parameterized by $\alpha$ with true value $\alpha_0$; that is,

$$W_i \equiv W(X_i; \alpha_0), \qquad i = 1, \ldots, n.$$

Let $\hat{\alpha}_n$ be an estimator of $\alpha$. Then, we can estimate $W_i$ by

$$\hat{W}_i = W(X_i; \hat{\alpha}_n), \qquad i = 1, \ldots, n.$$

In this subsection, we take $\Omega_i = W_i$, $i = 1, \ldots, n$, for simplicity, and we consider the asymptotic properties of the estimator $\hat{\theta}_n^*$, which are obtained from the following semiparametric estimating function with estimated weights:

$$(3.4) \quad \Psi_n^*(\theta, \hat{\eta}_n^*, \hat{\rho}_n^*) = \frac{1}{n} \sum_{i=1}^n \hat{W}_i \hat{\rho}_n^*(Y_i - \theta' Z_i, \theta) \{Z_i - \hat{\eta}_n^*(Y_i - \theta' Z_i, \theta)\} \Delta_i,$$

where

$$(3.5) \qquad \hat{\eta}_n^*(t, \theta) = \frac{\sum_{j=1}^n \hat{W}_j 1(Y_j - \theta' Z_j \geq t) Z_j}{\sum_{j=1}^n \hat{W}_j 1(Y_j - \theta' Z_j \geq t)}$$

and

$$(3.6) \qquad \hat{\rho}_n^*(t, \theta) = \frac{\sum_{j=1}^n \hat{W}_j 1(Y_j - \theta' Z_j \geq t)}{\sum_{j=1}^n \hat{W}_j}.$$

This case $\Omega_i = W_i$ handles the case-cohort study, naturally, when inverse sampling probability weights are used for which $\Omega_i = W_i = 1$ whenever $\Delta_i = 1$. Note that the estimating function (3.4) is obtained by replacing known weights $W_i$ with their estimates $\hat{W}_i$ in $\Psi_n(\theta, \hat{\eta}_n, \hat{\rho}_n)$, $\hat{\eta}_n$ and $\hat{\rho}_n$; see (2.2) and (2.3). As in Theorem 3.3, the following result holds for other types of bounded weight function $\rho_0$ and estimator $\hat{\rho}_n^*$, provided that $\{\hat{\rho}_n^*\}$ and $\{\rho_0\}$ are Donsker, and that $\hat{\rho}_n^*$, as a function of $\alpha$, is an asymptotically linear estimator that is twice continuously differentiable in $\alpha$ with the first-order derivative converging to an integrable limit at $\alpha_0$. The latter remark becomes clear in the proof of the next theorem.

We now consider consistency and asymptotic normality of $\hat{\theta}_n^*$ in Theorem 3.4 with a reasonable assumption about $\hat{\alpha}_n$ and a classical smoothness condition for $W(X;\alpha)$ in $\alpha$. The efficiency gain from using estimated weights becomes evident.

THEOREM 3.4. *Suppose that $W(X;\alpha)$ is twice differentiable, with respect to $\alpha$, in $\mathcal{A}_0 \times \mathcal{X}$ with continuous and bounded derivatives, where $\mathcal{A}_0$ is a neighborhood of $\alpha_0$ and $\mathcal{X}$ is the bounded sample space of the random variable $X$. Suppose that $\hat{\alpha}_n$ is an asymptotically efficient estimator of $\alpha$ with bounded influence function at $\alpha_0$. Let $\hat{\eta}_n^*$ and $\hat{\rho}_n^*$ be defined by (3.5) and (3.6), and let $\hat{\theta}_n^*$ be an approximate root satisfying the equation $\Psi_n^*(\hat{\theta}_n^*, \hat{\eta}_n^*(\cdot,\hat{\theta}_n^*), \hat{\rho}_n^*(\cdot,\hat{\theta}_n^*)) = o_{p^*}(n^{-1/2})$. Suppose that all the assumptions in Theorem 3.3 hold. Then, $\hat{\theta}_n^*$ is consistent, and $n^{1/2}(\hat{\theta}_n^* - \theta_0)$ is asymptotically normal with zero mean and the asymptotic variance*

$$(3.7) \qquad \Sigma_0 - \{\dot{\Psi}_\theta(\theta_0,\eta_0,\rho_0)\}^{-1} B V_0 B' \{\dot{\Psi}_\theta(\theta_0,\eta_0,\rho_0)\}^{-1},$$

*where $\Sigma_0$ is the asymptotic variance of $n^{1/2}(\hat{\theta}_n - \theta_0)$ determined by (3.3), $V_0$ is the asymptotic variance of $n^{1/2}(\hat{\alpha}_n - \alpha_0)$, and*

$$B = P[\rho_0(\varepsilon_0,\theta_0) A_2(\varepsilon_0,\theta_0) \Delta] - P[\rho_0(\varepsilon_0,\theta_0)\{Z - \eta_0(\varepsilon_0,\theta_0)\}(\dot{W}_\alpha(X;\alpha_0))' \Delta],$$

*with $\dot{W}_\alpha(X;\alpha)$ denoting the $\alpha$-derivative of $W(X;\alpha)$ and*

$$A_2(t,\theta_0) = \frac{1}{\rho_0(t,\theta_0)}[P\{1(\varepsilon_0 \geq t) Z (\dot{W}_\alpha(X;\alpha_0))'\}$$
$$- \eta_0(t,\theta_0) P\{(\dot{W}_\alpha(X;\alpha_0))' 1(\varepsilon_0 \geq t)\}].$$

Note that, if $\hat{\rho}_n^* = \hat{\rho}_n = 1$, then $\rho_0(t,\theta_0)$, in the above expression for $A_2$, should be replaced by $P\{1(\varepsilon_0 \geq t)\}$. The asymptotic efficiency of $\hat{\alpha}_n$ is one of three sufficient conditions for applying the result of [18] to obtain the above asymptotic normality of $\hat{\theta}_n^*$. When data are missing at random and inverse sampling probability weights are considered, the parameter $\alpha$ is adaptive to other parameters (see [1]) and its efficient estimator can be easily obtained, for example, by the maximum likelihood method. In sampling designs, a stratified approach is commonly used to improve efficiency. If the number of strata is finite, then the (independent Bernoulli) sampling probabilities within strata consist of the parameter $\alpha$, and the sampling fractions are the maximum likelihood estimates of $\alpha$.

The other two conditions of [18] are: (i) $n^{1/2}(\hat{\theta}_n - \theta_0)$ and $n^{1/2}(\hat{\alpha}_n - \alpha_0)$ are asymptotically jointly normal; and (ii) $n^{1/2}(\hat{\theta}_n^* - \theta_0)$ is asymptotically equivalent to $n^{1/2}(\hat{\theta}_n - \theta_0) + B n^{1/2}(\hat{\alpha}_n - \alpha_0)$. The former is determined by (3.3) in Theorem 3.3 and the fact that $\hat{\alpha}_n$ is an asymptotically linear estimator. The latter is established with a detailed proof in Section 6.

Consider a stratified case-cohort study. Suppose that all the censored subjects in a study cohort are divided into $S$ strata by the variable $Z_2 \in \{\zeta_1, \ldots, \zeta_S\}$. In a stratified case-cohort study, all of the failures are completely observed. For censored subjects, we denote the true sampling probabilities by $\alpha_{0s}$, $1 \leq s \leq S$. Suppose that there are $n_s$ subjects in stratum $s$, out of whom $n_s^*$ are selected into the subcohort by the independent Bernoulli sampling. We assume that, when $n \to \infty$, $n_s/n \to \gamma_s > 0$, $1 \leq s \leq S$. Instead of using the true sampling probabilities $\alpha_0 = (\alpha_{01}, \ldots, \alpha_{0S})'$ in the weight function $W$, we now replace each $\alpha_{0s}$ with the sampling fraction $\hat{\alpha}_{n,s} = n_s^*/n_s$, $1 \leq s \leq S$. We can then denote the sampling probability and its estimator of the $i$th subject as

$$\pi_i = \sum_{s=1}^{S} 1(Z_{2i} = \zeta_s)\alpha_{0s} \quad \text{and} \quad \hat{\pi}_i = \sum_{s=1}^{S} 1(Z_{2i} = \zeta_s)\hat{\alpha}_{n,s}.$$

We consider the inverse sampling probability weights

$$W(X_i; \hat{\alpha}_n) = \Delta_i + (1 - \Delta_i)\frac{1(i \in \mathcal{SC})}{\hat{\pi}_i}.$$

The second term in the expression for matrix $B$ in Theorem 3.4 becomes zero, since $\dot{W}_\alpha$ contains the factor $(1 - \Delta)$. The asymptotic variance of $\hat{\alpha}_n$ is $V_0 = \text{diag}\{\alpha_{01}(1 - \alpha_{01})/\gamma_1, \ldots, \alpha_{0S}(1 - \alpha_{0S})/\gamma_S\}$, which can be easily estimated from observed data.

## 4. Numerical results.

4.1. *Asymptotic efficiency comparison.* Considering the standard normal, standard logistic and standard extreme value error distributions in model (1.1), we evaluate asymptotic efficiency under a case-cohort setting to illustrate different extents of efficiency gain by using different weights. The one-dimensional covariate $Z$ is taken to follow a Bernoulli distribution with success probability 0.3 and $\theta_0 = 0$. Censoring time has a uniform distribution on $[a, b]$, where $a$ and $b$ are chosen to obtain 80% censoring proportion. Let $Z^*$ be a binary correlate of $Z$ with $\Pr(Z^* = 1|Z = 1) = 0.8$ and $\Pr(Z^* = 0|Z = 0) = 0.8$. The subcohort is a stratified subsample selected by independent Bernoulli sampling with selection probability $\pi(Z^*)$, chosen so that the two strata determined by $Z^*$ have the same expected number of subjects.

For each error distribution, we consider a $2^3$ factorial design with the following factors:

- logrank weights ($\hat{\rho}_n = 1$) and Gehan weights [see (2.3)];
- subject specific weight: predictable with $W_i = 1(i \in \mathcal{SC})/\pi_i$ and nonpredictable with $W_i = \Delta_i + (1 - \Delta_i)1(i \in \mathcal{SC})/\pi_i$;

- subject specific weights: true $W_i = W(X_i; \alpha_0)$ and estimated $\hat{W}_i = W(X_i; \hat{\alpha}_n)$.

The asymptotic variance of logrank weighted method for the full cohort is used as the benchmark, and we report the relative efficiency for each of the 8 scenarios with subcohort size fraction ranging from 1% to 100%. Results are given in Figures 1–3, where: (1) dark curves represent logrank weights, and gray curves represent Gehan weights; (2) solid curves represent predictable known weights, and dotted curves represent predictable estimated weights; and (3) dashed curves represent nonpredictable known weights, and dotted/dashed curves represent nonpredictable estimated weights.

We can see that using estimated weights $W(X_i; \hat{\alpha}_n)$ does not improve efficiency very much compared to using true weights $W(X_i; \alpha_0)$ for the settings considered. The efficiency gain from using the nonpredictable weights is substantial, especially for small to moderate sampling rates. An interesting feature is that when the subcohort size is relatively small, the Gehan weighted method performs much better than the logrank weighted method for all three error distributions, even though the result is opposite when subcohort size is close to the full cohort for both logistic and extreme value error distributions. We do not have an analytical explanation for this phenomenon, which seems to persist in other simulations as well. It seems safe, however, to recommend the Gehan weights for the problems with missing data; it is fortuitous that the Gehan weights also yield a monotone estimating function, which is a numerically advantageous property. Another interesting phenomenon is that, for the logistic error, the Gehan weights may be somewhat less efficient than the logrank weights for censored data, even though they are the most efficient for uncensored data (see [12]).

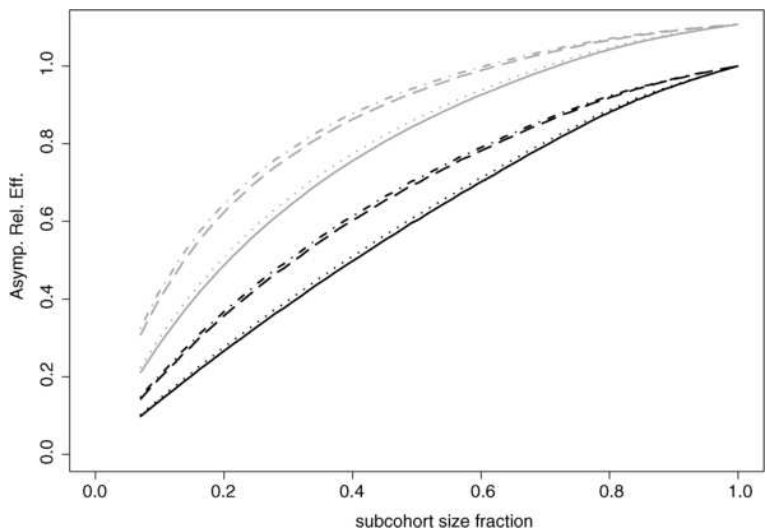

FIG. 1. *Asymptotic efficiency under normal error distribution.*

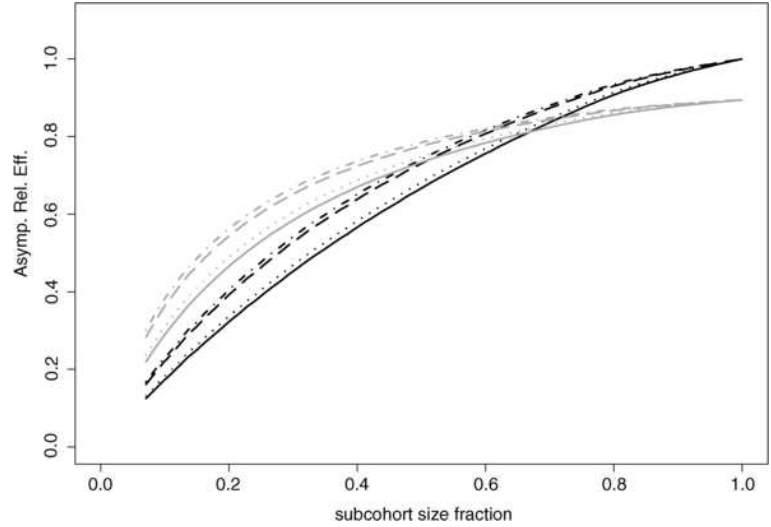

Fig. 2. *Asymptotic efficiency under logistic error distribution.*

4.2. *Simulations.* We conduct simulations under the same settings as that in the previous subsection. Since the simulation results are basically telling the same story for different error distributions, we only report the re-sults for the logistic error. We consider case-cohort designs with cohort size of 2000 and subcohort sizes of 15%, 20% and 25% of the entire cohort on average, which lead to on average 640, 720 and 800 completely observed subjects, respectively. Bias of the point estimator, average of the variance estimator, empirical variance and 95% coverage probability, based on the variance estimator, are reported for five different analyses using the following logrank and Gehan weights: full data analysis, predictable subject-specific weighted analysis using true weights, predictable subject-specific weighted analysis using estimated weights, nonpredictable subject-specific weighted analysis using true weights and nonpredictable subject-specific weighted analysis using estimated weights. The asymptotic variance for each scenario is also reported. From Table 1, we see that all of the methods work well for finite samples and reflect the patterns observed from the efficiency results in the previous subsection.

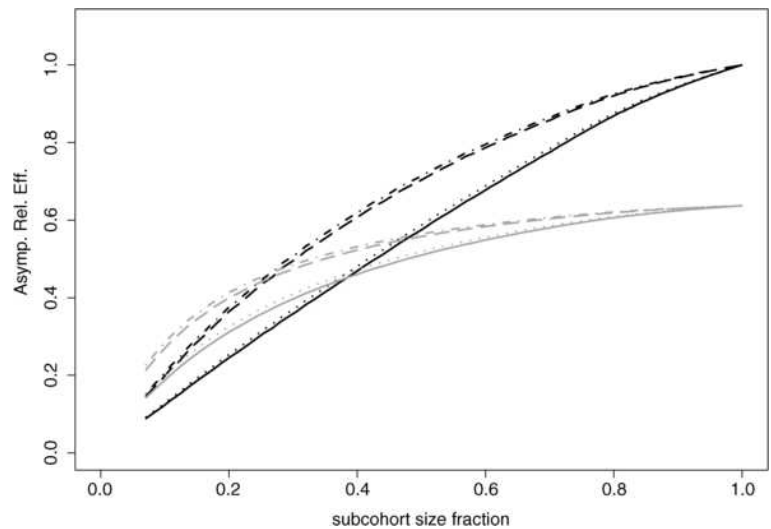

Fig. 3. *Asymptotic efficiency under extreme value error distribution.*

**5. Discussion.** We consider only the case where weights $\Omega_i$ and $W_i$ are i.i.d. for all $i = 1, \ldots, n$, which makes the proofs of the asymptotic properties more straightforward. For the case where the weights are determined by (stratified) simple random sampling, the method of [3] may be applicable, and this is an interesting topic worthy of further investigation.

**6. Proofs.**

6.1. *Proof of Theorem 3.1.* As in [26], for notational simplicity, we assume one-dimensional $\theta$ in the proofs of the theorems in Section 3.

Since $\eta$, $\eta_0$, $\rho$ and $\rho_0$ belong to Glivenko–Cantelli classes, it follows, from Theorem 3 of [28], that the set of bounded functions $\{\Omega\rho(Y,\theta)\{Z - \eta(\varepsilon_\theta,\theta)\}\Delta\}$ is a Glivenko–Cantelli class. By adding and subtracting the same term, and by the triangle inequality, we then have that

$$\|\Psi_n(\theta,\eta,\rho) - \Psi(\theta,\eta_0,\rho_0)\|$$
$$= \|\mathbb{P}_n[\Omega\rho(\varepsilon_\theta,\theta)\{Z - \eta(\varepsilon_\theta,\theta)\}\Delta] - P[\Omega\rho_0(\varepsilon_\theta,\theta)\{Z - \eta_0(\varepsilon_\theta,\theta)\}\Delta]\|$$
$$\leq \|(\mathbb{P}_n - P)[\Omega\rho(\varepsilon_\theta,\theta)\{Z - \eta(\varepsilon_\theta,\theta)\}\Delta]\|$$
$$+ \|P\{\Omega(\rho - \rho_0)Z\Delta\}\| + \|P\{\Omega(\rho\eta - \rho_0\eta_0)\Delta\}\|.$$

The first term on the right-hand side of the above inequality converges to zero in outer probability by the Glivenko–Cantelli property. Obviously,

$$\|P\{\Omega(\rho - \rho_0)Z\Delta\}\| \leq \|\rho - \rho_0\|P|\Omega Z\Delta| \to 0$$

and

$$\|P\{\Omega(\rho\eta - \rho_0\eta_0)\Delta\}\|$$
$$\leq \|\rho\eta - \rho_0\eta_0\|P|\Omega\Delta| \to 0,$$

TABLE 1
*Summary statistics of simulations, where $\alpha$ = subcohort size fraction; Method 1 = full data analysis, 2 = predictable subject-specific weighted analysis using true weights, 3 = predictable subject-specific weighted analysis using estimated weights, 4 = nonpredictable subject-specific weighted analysis using true weights, 5 = nonpredictable subject-specific weighted analysis using estimated weights; Emp. Var = empirical variance estimator; Ave. Var = average of variance estimator; CP = coverage probability; Asym. Var = asymptotic variance*

| $\alpha$ | Weight | Method | $\hat{\theta}_n$ | Emp. Var | Ave. Var | 95% CP | Asym. Var |
|---|---|---|---|---|---|---|---|
| 0.15 | Logrank | 1 | −0.001 | 0.018 | 0.019 | 95.6 | 0.018 |
|  |  | 2 | 0.019 | 0.074 | 0.075 | 93.2 | 0.073 |
|  |  | 3 | 0.018 | 0.066 | 0.072 | 94.4 | 0.069 |
|  |  | 4 | 0.015 | 0.056 | 0.059 | 95.4 | 0.059 |
|  |  | 5 | 0.015 | 0.052 | 0.058 | 95.8 | 0.056 |
|  | Gehan | 1 | 0.006 | 0.020 | 0.020 | 96.6 | 0.020 |
|  |  | 2 | 0.018 | 0.047 | 0.047 | 94.0 | 0.047 |
|  |  | 3 | 0.016 | 0.040 | 0.042 | 95.4 | 0.044 |
|  |  | 4 | 0.015 | 0.038 | 0.039 | 96.4 | 0.039 |
|  |  | 5 | 0.014 | 0.034 | 0.036 | 96.2 | 0.037 |
| 0.20 | Logrank | 1 | −0.001 | 0.018 | 0.019 | 95.6 | 0.018 |
|  |  | 2 | 0.007 | 0.060 | 0.059 | 94.0 | 0.056 |
|  |  | 3 | 0.008 | 0.055 | 0.057 | 94.8 | 0.054 |
|  |  | 4 | 0.006 | 0.049 | 0.048 | 93.0 | 0.046 |
|  |  | 5 | 0.007 | 0.046 | 0.046 | 94.6 | 0.045 |
|  | Gehan | 1 | 0.006 | 0.020 | 0.020 | 96.6 | 0.020 |
|  |  | 2 | 0.011 | 0.039 | 0.039 | 96.0 | 0.039 |
|  |  | 3 | 0.012 | 0.035 | 0.035 | 95.6 | 0.037 |
|  |  | 4 | 0.011 | 0.034 | 0.033 | 95.2 | 0.033 |
|  |  | 5 | 0.011 | 0.031 | 0.031 | 95.8 | 0.032 |
| 0.25 | Logrank | 1 | −0.001 | 0.018 | 0.019 | 95.6 | 0.018 |
|  |  | 2 | 0.003 | 0.048 | 0.049 | 94.0 | 0.047 |
|  |  | 3 | 0.004 | 0.043 | 0.047 | 95.8 | 0.045 |
|  |  | 4 | 0.002 | 0.040 | 0.041 | 94.4 | 0.039 |
|  |  | 5 | 0.003 | 0.038 | 0.040 | 94.8 | 0.038 |
|  | Gehan | 1 | 0.006 | 0.020 | 0.020 | 96.6 | 0.020 |
|  |  | 2 | 0.007 | 0.034 | 0.034 | 95.6 | 0.034 |
|  |  | 3 | 0.008 | 0.031 | 0.032 | 95.4 | 0.033 |
|  |  | 4 | 0.008 | 0.030 | 0.030 | 95.0 | 0.030 |
|  |  | 5 | 0.008 | 0.028 | 0.029 | 94.8 | 0.029 |

where

$$\|\rho\eta - \rho_0\eta_0\| = \tfrac{1}{2}\|(\rho - \rho_0)(\eta + \eta_0) + (\rho + \rho_0)(\eta - \eta_0)\|$$
$$\leq \tfrac{1}{2}\|\rho - \rho_0\| \cdot \|\eta + \eta_0\| + \tfrac{1}{2}\|\rho + \rho_0\| \cdot \|\eta - \eta_0\|$$

$$\to 0.$$

This establishes (3.1), which, in turn, can be shown to imply $|\hat{\theta}_n - \theta_0| \to 0$ in outer probability, as in [8]. For completeness, we include the argument here.

Since $\theta_0$ is the unique solution to $\Psi(\theta, \eta_0(\cdot, \theta), \rho_0(\cdot, \theta)) = 0$, for any fixed $\varepsilon > 0$, there exists a $\delta > 0$ such that

$$P[|\hat{\theta}_n - \theta_0| > \varepsilon] \leq P[|\Psi(\hat{\theta}_n, \eta_0(\cdot, \hat{\theta}_n), \rho_0(\cdot, \hat{\theta}_n))| > \delta].$$

We show that $|\Psi(\hat{\theta}_n, \eta_0(\cdot, \hat{\theta}_n), \rho_0(\cdot, \hat{\theta}_n))| \to 0$ in outer probability, and the consistency of $\hat{\theta}_n$ follows immediately. Note that there exists a sequence $\{\delta_n\} \downarrow 0$ such that $\|\eta - \eta_0\| \leq \delta_n$ and $\|\rho - \rho_0\| \leq \delta_n$ with probability tending to one. Hence, from (3.1), we have the inequalities

$$\begin{aligned}
|\Psi(\hat{\theta}_n, \eta_0(\cdot, \hat{\theta}_n)), \rho_0(\cdot, \hat{\theta}_n))| \\
\leq |\Psi_n(\hat{\theta}_n, \eta(\cdot, \hat{\theta}_n), \rho(\cdot, \hat{\theta}_n))| \\
+ |\Psi(\hat{\theta}_n, \eta_0(\cdot, \hat{\theta}_n), \rho_0(\cdot, \hat{\theta}_n)) - \Psi_n(\hat{\theta}_n, \eta(\cdot, \hat{\theta}_n), \rho(\cdot, \hat{\theta}_n))| \\
\leq |\Psi_n(\hat{\theta}_n, \eta(\cdot, \hat{\theta}_n), \rho(\cdot, \hat{\theta}_n))| + o_{p^*}(1) \\
= o_{p^*}(1).
\end{aligned}$$

Hence, $\hat{\theta}_n$ is consistent.

We now show that (3.1) holds, when $\eta$ and $\rho$ are replaced by $\hat{\eta}_n$ and $\hat{\rho}_n$ given in (2.2) and (2.3), respectively, and $\rho_0(t, \theta) = \Pr(\varepsilon_\theta \geq t)$. We define

$$D_n^{(0)}(t, \theta) = \mathbb{P}_n\{W1(\varepsilon_\theta \geq t)\}, \qquad d^{(0)}(t, \theta) = P\{W1(\varepsilon_\theta \geq t)\};$$
$$D_n^{(1)}(t, \theta) = \mathbb{P}_n\{W1(\varepsilon_\theta \geq t)Z\}, \qquad d^{(1)}(t, \theta) = P\{W1(\varepsilon_\theta \geq t)Z\}.$$

Thus, $\hat{\eta}_n(t, \theta) = D_n^{(1)}(t, \theta)/D_n^{(0)}(t, \theta)$ and $\eta_0(t, \theta) = d^{(1)}(t, \theta)/d^{(0)}(t, \theta)$. The latter equality holds because

$$P\{W1(\varepsilon_\theta \geq t)\} = P\{1(\varepsilon_\theta \geq t)\} \quad \text{and} \quad P\{W1(\varepsilon_\theta \geq t)Z\} = P\{1(\varepsilon_\theta \geq t)Z\}.$$

Since the class of functions $\{1(\varepsilon_\theta \geq t)\}$ is a VC-class (see, e.g., Exercise 9 on page 151 and Exercise 14 on page 152 in [27]) and, thus, a Donsker class, we know that the sets of functions $\mathcal{F}_0 = \{W1(\varepsilon_\theta \geq t)\}$ and $\mathcal{F}_1 = \{W1(\varepsilon_\theta \geq t)Z\}$ are Donsker classes (see, e.g., [27], Section 2.10). Since Donsker classes are Glivenko–Cantelli classes, it follows that $\|D_n^{(k)}(t, \theta) - d^{(k)}(t, \theta)\| \to 0$ in outer probability, $k = 0, 1$. Let $\tau$ correspond to $T^*$ in [26] and represent the longest follow-up time. Since both $D_n^{(0)}$ (with probability 1) and $d^{(0)}$ are bounded away from zero when $t \leq \tau$, we have

$$\|\hat{\eta}_n(t, \theta) - \eta_0(t, \theta)\| \to 0 \tag{6.1}$$

in outer probability. Similarly, we have

(6.2) $$\|\hat{\rho}_n(t,\theta) - \rho_0(t,\theta)| \to 0$$

in outer probability.

Let $\bar{\mathcal{F}}_k$ be the closure of $\mathcal{F}_k$, $k = 0, 1$, respectively, in which the convergence is both pointwise and in $L_2(P)$. Then, $D_n^{(k)}(t,\theta)$ and $d^{(k)}(t,\theta)$ are in the convex hull of $\bar{\mathcal{F}}_k$, $k = 0, 1$, and, thus, belong to Donsker classes (see, e.g., [27], Theorems 2.10.2 and 2.10.3). Hence, both $\{\hat{\eta}_n(t,\theta)\}$ and $\{\eta_0(t,\theta)\}$ are Donsker (by [27], Example 2.10.9) and, thus, Glivenko–Cantelli. Similarly, we can argue that both $\{\hat{\rho}_n(t,\theta)\}$ and $\{\rho_0(t,\theta)\}$ are Donsker and, hence, Glivenko–Cantelli. Then, by the first half of the proof we obtain

$$\|\Psi_n(\theta, \hat{\eta}_n, \hat{\rho}_n) - \Psi(\theta, \eta_0, \rho_0)\| \to 0$$

in outer probability.

6.2. *Proof of Theorem 3.2.* From the proof of Theorem 3.1 we see that $n^{1/2}\{D_n^{(k)}(t,\theta) - d^{(k)}(t,\theta)\}$, $k = 0, 1$, converge to zero mean Gaussian processes for all $\theta \in \Theta_0$, and $\|n^{1/2}\{D_n^{(k)}(t,\theta) - d^{(k)}(t,\theta)\}\|_0 = O_{p^*}(1)$, $k = 0, 1$, by the tail bounds for the supremum of empirical processes in [27], Section 2.14. We then have

$$\begin{aligned}
&n^{1/2}\{\hat{\eta}_n(t,\theta) - \eta_0(t,\theta)\} \\
&= n^{1/2}\bigg[\frac{1}{d^{(0)}(t,\theta)}\{D_n^{(1)}(t,\theta) - d^{(1)}(t,\theta)\} \\
&\qquad - \frac{D_n^{(1)}(t,\theta)}{D_n^{(0)}(t,\theta)d^{(0)}(t,\theta)}\{D_n^{(0)}(t,\theta) - d^{(0)}(t,\theta)\}\bigg] \\
&= n^{1/2}\bigg[\frac{1}{d^{(0)}(t,\theta)}\{D_n^{(1)}(t,\theta) - d^{(1)}(t,\theta)\} \\
&\qquad - \frac{d^{(1)}(t,\theta)}{d^{(0)}(t,\theta)^2}\{D_n^{(0)}(t,\theta) - d^{(0)}(t,\theta)\}\bigg] + o_{p^*}(1) \\
&= d^{(0)}(t,\theta)^{-1}n^{1/2}[\{D_n^{(1)}(t,\theta) - D_n^{(0)}(t,\theta)\eta_0(t,\theta)\} \\
&\qquad - \{d^{(1)}(t,\theta) - d^{(0)}(t,\theta)\eta_0(t,\theta)\}] + o_{p^*}(1) \\
&= d^{(0)}(t,\theta)^{-1}\mathbb{G}_n[W1(\varepsilon_\theta \geq t)\{Z - \eta_0(t,\theta)\}] + o_{p^*}(1).
\end{aligned}$$

Since the classes of functions $\{W\}$, $\{1(\varepsilon_\theta \geq t)\}$, $\{Z\}$ and $\{\eta_0\}$ are Donsker, we know that $\{W1(\varepsilon_\theta \geq t)\{Z - \eta_0(t,\theta)\}\}$ is Donsker (e.g., [27], Section 2.10). Thus, $n^{1/2}\|\hat{\eta}_n - \eta_0\|_0 = O_{p^*}(1)$, since $d^{(0)}(t,\theta)^{-1}$ is bounded.

We now show $n^{1/2}|\hat{\theta}_n - \theta| = O_{p^*}(1)$. First, we have

(6.3) $$\|n^{1/2}\{\Psi_n(\theta, \hat{\eta}_n(\cdot,\theta), \hat{\rho}_n(\cdot,\theta)) - \Psi(\theta, \eta_0(\cdot,\theta), \rho_0(\cdot,\theta))\}\|_0 = O_{p^*}(1)$$

by applying the triangle inequality, and that $\{\hat{\eta}_n\}$ and $\{\hat{\rho}_n\}$ are Donsker, as well as $n^{1/2}\|\hat{\rho}_n - \rho_0\|_0 = O_{p^*}(1)$ and $n^{1/2}\|\hat{\eta}_n - \eta_0\|_0 = O_{p^*}(1)$ in the following calculation:

$$\|n^{1/2}\{\Psi_n(\theta, \hat{\eta}_n(\cdot, \theta), \hat{\rho}_n(\cdot, \theta)) - \Psi(\theta, \eta_0(\cdot, \theta), \rho_0(\cdot, \theta))\}\|_0$$
$$= \|n^{1/2}(\mathbb{P}_n - P)[\Omega\hat{\rho}_n(\varepsilon_\theta, \theta)\{Z - \hat{\eta}_n(\varepsilon_\theta, \theta)\}\Delta]$$
$$+ n^{1/2}P[\Omega\{\hat{\rho}_n(\varepsilon_\theta, \theta) - \rho_0(\varepsilon_\theta, \theta)\}Z\Delta]$$
$$+ n^{1/2}P[\Omega\hat{\rho}_n(\varepsilon_\theta, \theta)\hat{\eta}_n(\varepsilon_\theta, \theta) - \rho_0(\varepsilon_\theta, \theta)\eta_0(\varepsilon_\theta, \theta)\Delta]\|_0$$
$$\leq \|\mathbb{G}_n[\Omega\hat{\rho}_n(\varepsilon_\theta, \theta)\{Z - \hat{\eta}_n(\varepsilon_\theta, \theta)\}\Delta]\|_0 + n^{1/2}\|\hat{\rho}_n - \rho_0\|_0 \cdot P(\Omega Z\Delta)$$
$$+ \tfrac{1}{2}(n^{1/2}\|\hat{\rho}_n - \rho_0\|_0 \cdot \|\hat{\eta}_n + \eta_0\|_0$$
$$+ \|\hat{\rho}_n + \rho_0\|_0 \cdot n^{1/2}\|\hat{\eta}_n - \eta_0\|_0)P(\Omega\Delta)$$
$$= O_{p^*}(1).$$

Because $\Psi(\theta_0, \eta_0(\cdot, \theta_0), \rho_0(\cdot, \theta_0)) = 0$ and $|\hat{\theta}_n - \theta_0| = o_{p^*}(1)$ by Theorem 3.1, we then have

$$O_{p^*}(1) = -n^{1/2}\{\Psi_n(\hat{\theta}_n, \hat{\eta}_n(\cdot, \hat{\theta}_n), \hat{\rho}_n(\cdot, \hat{\theta}_n)) - \Psi(\hat{\theta}_n, \eta_0(\cdot, \hat{\theta}_n), \rho_0(\cdot, \hat{\theta}_n))\}$$
$$= o_{p^*}(1) + n^{1/2}\Psi(\hat{\theta}_n, \eta_0(\cdot, \hat{\theta}_n), \rho_0(\cdot, \hat{\theta}_n))$$
(6.4)
$$\quad - n^{1/2}\Psi(\theta_0, \eta_0(\cdot, \theta_0), \rho_0(\cdot, \theta_0))$$
$$= o_{p^*}(1) + n^{1/2}(\hat{\theta}_n - \theta_0)\dot{\Psi}_\theta(\theta^*, \eta_0(\cdot, \theta^*), \rho_0(\cdot, \theta^*))$$
$$= o_{p^*}(1) + n^{1/2}(\hat{\theta}_n - \theta_0)\{\dot{\Psi}_\theta(\theta_0, \eta_0(\cdot, \theta_0), \rho_0(\cdot, \theta_0)) + o_{p^*}(1)\},$$

where $\theta^*$ is a point between $\theta_0$ and $\hat{\theta}_n$. Thus, $n^{1/2}(\hat{\theta}_n - \theta_0) = O_{p^*}(1)$.

Let $C_n = n^{-1}\sum_{i=1}^n W_i$. By the central limit theorem, $n^{1/2}(C_n - 1) = O_p(1)$. Thus, when $\hat{\rho}_n$ takes the form, in (2.3) and $\rho_0(t, \theta) = \Pr(\varepsilon_\theta \geq t)$, they are clearly bounded, and we can show $n^{1/2}\|\hat{\rho}_n - \rho_0\|_0 = O_{p^*}(1)$ by the following calculation:

$$n^{1/2}\{\hat{\rho}_n(t, \theta) - \rho_0(t, \theta)\}$$
$$= n^{1/2}\left[\{D_n^{(0)}(t, \theta) - d^{(0)}(t, \theta)\} - \frac{D_n^{(0)}(t, \theta)}{C_n}\{C_n - 1\}\right]$$
$$= n^{1/2}[\{D_n^{(0)}(t, \theta) - d^{(0)}(t, \theta)\} - d^{(0)}(t, \theta)\{C_n - 1\}] + o_{p^*}(1)$$
$$= n^{1/2}[\{D_n^{(0)}(t, \theta) - C_n d^{(0)}(t, \theta)\}] + o_{p^*}(1)$$
$$= \mathbb{G}_n[W\{1(\varepsilon_\theta \geq t) - d^{(0)}(t, \theta)\}] + o_{p^*}(1).$$

We have already shown in the proof of Theorem 3.1 that such chosen $\hat{\rho}_n$ and $\rho_0$ belong to a Donsker class.

6.3. *Proof of Theorem 3.3.* The differentiability of both $\rho_0(\varepsilon_\theta, \theta)$ and $\eta_0(\varepsilon_\theta, \theta)$ in $\theta$ and its implication of the differentiability of $\Psi(\theta, \eta_0(\cdot, \theta), \rho_0(\cdot, \theta))$ in $\theta$, as well as the continuity and boundedness of the derivatives, can be shown by interchanging integration and differentiation, which is warranted by the dominated convergence theorem under the given regularity conditions. From Theorem 3.2, we know that $|\hat{\theta}_n - \theta_0| = O_{p^*}(n^{-1/2})$. Let $|\theta - \theta_0| \leq Kn^{-1/2}$ with $K < \infty$. Then, we have

$$
\begin{aligned}
& n^{1/2}\{\Psi_n(\theta, \hat{\eta}_n(\cdot, \theta), \hat{\rho}_n(\cdot, \theta)) - \Psi_n(\theta_0, \hat{\eta}_n(\cdot, \theta_0), \hat{\rho}_n(\cdot, \theta_0))\} \\
& \quad = n^{1/2}[\mathbb{P}_n\Omega\hat{\rho}_n(\varepsilon_\theta, \theta)\{Z - \hat{\eta}_n(\varepsilon_\theta, \theta)\}\Delta \quad (6.5) \\
& \qquad - \mathbb{P}_n\Omega\hat{\rho}_n(\varepsilon_\theta, \theta)\{Z - \hat{\eta}_n(\varepsilon_0, \theta_0)\}\Delta] \\
& \quad + n^{1/2}[\mathbb{P}_n\Omega\hat{\rho}_n(\varepsilon_\theta, \theta)\{Z - \hat{\eta}_n(\varepsilon_0, \theta_0)\}\Delta \quad (6.6) \\
& \qquad - \mathbb{P}_n\Omega\hat{\rho}_n(\varepsilon_0, \theta_0)\{Z - \hat{\eta}_n(\varepsilon_0, \theta_0)\}\Delta].
\end{aligned}
$$

We first look at term (6.5), which can be rewritten as

$$
\begin{aligned}
& n^{1/2}[-\mathbb{P}_n\Omega\hat{\rho}_n(\varepsilon_\theta, \theta)\hat{\eta}_n(\varepsilon_\theta, \theta)\Delta + \mathbb{P}_n\Omega\hat{\rho}_n(\varepsilon_\theta, \theta)\hat{\eta}_n(\varepsilon_0, \theta_0)\Delta] \\
(6.7) \quad & = -\mathbb{G}_n[\Omega\hat{\rho}_n(\varepsilon_\theta, \theta)\{\hat{\eta}_n(\varepsilon_\theta, \theta) - \hat{\eta}_n(\varepsilon_0, \theta_0)\}\Delta] \\
(6.8) \quad & \quad - n^{1/2}P[\Omega\hat{\rho}_n(\varepsilon_\theta, \theta)\{\hat{\eta}_n(\varepsilon_\theta, \theta) - \hat{\eta}_n(\varepsilon_0, \theta_0)\}\Delta].
\end{aligned}
$$

Term (6.7) converges to zero in outer probability, because $\Omega\hat{\rho}_n\hat{\eta}_n\Delta$ belongs to a Donsker class by arguments similar to those in the proof of Theorem 3.1, and $\Omega\hat{\rho}_n(\varepsilon_\theta, \theta)\{\hat{\eta}_n(\varepsilon_\theta, \theta) - \hat{\eta}_n(\varepsilon_0, \theta_0)\}\Delta$ converges to zero in quadratic mean. Let $t' = t - (\theta - \theta_0)z$. Direct calculation yields

$$
\begin{aligned}
& n^{1/2}P[\Omega\hat{\rho}_n(\varepsilon_\theta, \theta)\{\hat{\eta}_n(\varepsilon_\theta, \theta) - \eta_0(\varepsilon_\theta, \theta)\}\Delta] \\
& = n^{1/2}P\bigg[\hat{\rho}_n(\varepsilon_\theta, \theta)\bigg\{\frac{D_n^{(1)}(\varepsilon_\theta, \theta)}{D_n^{(0)}(\varepsilon_\theta, \theta)} - \frac{d^{(1)}(\varepsilon_\theta, \theta)}{d^{(0)}(\varepsilon_\theta, \theta)}\bigg\}\Delta\bigg] \\
& = n^{1/2}\int \hat{\rho}_n(t', \theta)\bigg[\frac{1}{d^{(0)}(t', \theta)}\{D_n^{(1)}(t', \theta) - d^{(1)}(t', \theta)\} \\
& \qquad\qquad\qquad \times \frac{D_n^{(1)}(t', \theta)}{D_n^{(0)}(t', \theta)d^{(0)}(t', \theta)}\{D_n^{(0)}(t', \theta) - d^{(0)}(t', \theta)\}\bigg] \\
(6.9) \quad & \qquad \times \delta\, dP_{\varepsilon_0, \Delta, Z}(t, \delta, z) \\
& = n^{1/2}\int \hat{\rho}_n(t', \theta)\bigg[\frac{1}{d^{(0)}(t', \theta)}\{D_n^{(1)}(t', \theta) - d^{(1)}(t', \theta)\} \\
& \qquad\qquad\qquad \times \frac{d^{(1)}(t', \theta)}{d^{(0)}(t', \theta)^2}\{D_n^{(0)}(t', \theta) - d^{(0)}(t', \theta)\}\bigg]
\end{aligned}
$$

$$\times \delta \, dP_{\varepsilon_0,\Delta,Z}(t,\delta,z) + o_{p^*}(1)$$
$$= \int \mathbb{G}_n \hat{\rho}_n(t',\theta) d^{(0)}(t',\theta)^{-1} W 1(\varepsilon_\theta \geq t')$$
$$\times \{Z - \eta_0(t',\theta)\} \, dP_{\varepsilon_0,\Delta,Z}(t,1,z) + o_{p^*}(1)$$
$$= \int \mathbb{G}_n \hat{\rho}_n(t',\theta) \ell(t',\theta,W,Z,\varepsilon_\theta) \, dP_{\varepsilon_0,\Delta,Z}(t,1,z) + o_{p^*}(1)$$

where $\ell(t',\theta,W,Z,\varepsilon_\theta) = d^{(0)}(t',\theta)^{-1} W 1(\varepsilon_\theta \geq t') \{Z - \eta_0(t',\theta)\}$ and $P_{\varepsilon_0,\Delta,Z}$ denotes the joint probability law of $(\varepsilon_0,\Delta,Z)$. Clearly, the class of functions $\{\hat{\rho}_n(t,\theta)\ell(t,\theta,W,Z,\varepsilon_\theta)\}$ is Donsker. The above middle equality holds because

$$\left| n^{1/2} \int \hat{\rho}_n(t',\theta) \left[ \frac{1}{d^{(0)}(t',\theta)} \{D_n^{(1)}(t',\theta) - d^{(1)}(t',\theta)\} \right. \right.$$
$$\left. - \frac{D_n^{(1)}(t',\theta)}{D_n^{(0)}(t',\theta) \, d^{(0)}(t',\theta)} \{D_n^{(0)}(t',\theta) - d^{(0)}(t',\theta)\} \right]$$
$$\times \delta \, dP_{\varepsilon_0,\Delta,Z}(t,\delta,z)$$
$$- n^{1/2} \int \hat{\rho}_n(t',\theta) \left[ \frac{1}{d^{(0)}(t',\theta)} \{D_n^{(1)}(t',\theta) - d^{(1)}(t',\theta)\} \right.$$
$$\left. - \frac{d^{(1)}(t',\theta)}{d^{(0)}(t',\theta)^2} \{D_n^{(0)}(t',\theta) - d^{(0)}(t',\theta)\} \right]$$
$$\left. \times \delta \, dP_{\varepsilon_0,\Delta,Z}(t,\delta,z) \right|$$
$$= \left| \int \hat{\rho}_n(t',\theta) \left\{ \frac{d^{(1)}(t',\theta)}{d^{(0)}(t',\theta)^2} - \frac{D_n^{(1)}(t',\theta)}{D_n^{(0)}(t',\theta) d^{(0)}(t',\theta)} \right\} \right.$$
$$\left. \times n^{1/2} \{D_n^{(0)}(t',\theta) - d^{(0)}(t',\theta)\} \delta \, dP_{\varepsilon_0,\Delta,Z}(t,\delta,z) \right|$$
$$\leq 1 \cdot \left\| \frac{d^{(1)}(t,\theta)}{d^{(0)}(t,\theta)^2} - \frac{D_n^{(1)}(t,\theta)}{D_n^{(0)}(t,\theta) d^{(0)}(t,\theta)} \right\|$$
$$\cdot \| n^{1/2} \{D_n^{(0)}(t,\theta) - d^{(0)}(t,\theta)\} \| \cdot 1$$
$$= o_{p^*}(1) \cdot O_{p^*}(1) \cdot 1 = o_{p^*}(1)$$

by the tail bounds for the supremum of empirical processes in [27], Section 2.14. Similarly, we have
$$n^{1/2} P[\Omega \hat{\rho}_n(\varepsilon_\theta,\theta)\{\hat{\eta}_n(\varepsilon_0,\theta_0) - \eta_0(\varepsilon_0,\theta_0)\}\Delta]$$
$$= \int \mathbb{G}_n \hat{\rho}_n(t',\theta) \ell(t,\theta_0,W,Z,\varepsilon_0) \, dP_{\varepsilon_0,\Delta,Z}(t,1,z) + o_{p^*}(1).$$

Thus, (6.8) becomes

$$
\begin{aligned}
&- n^{1/2}P[\hat{\rho}_n(\varepsilon_\theta,\theta)\{\eta_0(\varepsilon_\theta,\theta) - \eta_0(\varepsilon_0,\theta_0)\}\Delta] \\
&\quad + \int \mathbb{G}_n\hat{\rho}_n(t',\theta) \\
&\qquad \times \{\ell(t',\theta,W,Z,\varepsilon_\theta) - \ell(t,\theta_0,W,Z,\varepsilon_0)\}\,dP_{\varepsilon_0,\Delta,Z}(t,1,z) \\
&\quad + o_{p^*}(1).
\end{aligned}
\tag{6.10}
$$

Note that $n^{1/2}\{\eta_0(\varepsilon_\theta,\theta) - \eta_0(\varepsilon_0,\theta_0)\} = n^{1/2}(\theta - \theta_0)\{\dot\eta_{0\theta}(\varepsilon_{\theta^*},\theta^*)\}$ is bounded (by assumptions of bounded density functions for failure and censoring times in [30]), where $\dot\eta_{0\theta}$ denotes the derivative of $\eta_0$ with respect to $\theta$, and $\theta^*$ is a point between $\theta_0$ and $\theta$. Thus, by repeatedly using the dominate convergence theorem, we know that the first term in (6.10) equals

$$-n^{1/2}(\theta-\theta_0)P\{\rho_0(\varepsilon_\theta,\theta)\dot\eta_{0\theta}(\varepsilon_0,\theta_0)\Delta\} + o_{p^*}(1),$$

which in turn equals

$$-n^{1/2}(\theta-\theta_0)P\{\rho_0(\varepsilon_0,\theta_0)\dot\eta_{0\theta}(\varepsilon_0,\theta_0)\Delta\} + o_{p^*}(1).$$

It can be verified that $\hat{\rho}_n(t',\theta)\{\ell(t',\theta,W,Z,\varepsilon_\theta) - \ell(t,\theta_0,W,Z,\varepsilon_0)\}$ converges to zero in quadratic mean; thus,

$$\|\mathbb{G}_n\hat{\rho}_n(t',\theta)\{\ell(t',\theta,W,Z,\varepsilon_\theta) - \ell(t,\theta_0,W,Z,\varepsilon_0)\}\| = o_{p^*}(1),$$

then the second term in (6.10) converges to zero in outer probability. So we have shown that term (6.5) is asymptotically equivalent to $-n^{1/2}(\theta-\theta_0) \times P\{\rho_0(\varepsilon_0,\theta_0)\dot\eta_{0\theta}(\varepsilon_0,\theta_0)\Delta\}$.

We now consider term (6.6), which can be rewritten as

$$n^{1/2}\mathbb{P}_n[\Omega\{Z-\hat\eta_n(\varepsilon_0,\theta_0)\}\Delta\{\hat\rho_n(\varepsilon_\theta,\theta) - \hat\rho_n(\varepsilon_0,\theta_0)\}]$$

$$= \mathbb{G}_n[\Omega\{Z-\hat\eta_n(\varepsilon_0,\theta_0)\}\Delta\{\hat\rho_n(\varepsilon_\theta,\theta) - \hat\rho_n(\varepsilon_0,\theta_0)\}] \tag{6.11}$$

$$\quad + n^{1/2}P[\Omega\{Z-\hat\eta_n(\varepsilon_0,\theta_0)\}\Delta\{\hat\rho_n(\varepsilon_\theta,\theta) - \hat\rho_n(\varepsilon_0,\theta_0)\}]. \tag{6.12}$$

Because $\Omega\{Z-\hat\eta_n(\varepsilon_0,\theta_0)\}\Delta\{\hat\rho_n(\varepsilon_\theta,\theta) - \hat\rho_n(\varepsilon_0,\theta_0)\}$ belongs to a Donsker class and converges to zero in quadratic mean, we know that term (6.11) converges to zero in outer probability. Similar to the calculation in (6.9), for (6.12), we have

$$
\begin{aligned}
&n^{1/2}P[\Omega\{Z-\hat\eta_n(\varepsilon_0,\theta_0)\}\Delta\{\hat\rho_n(\varepsilon_\theta,\theta) - \rho_0(\varepsilon_\theta,\theta)\}] \\
&= n^{1/2}\int \{z-\hat\eta_n(t,\theta_0)\}\bigg[\{D_n^{(0)}(t',\theta) - d^{(0)}(t',\theta)\} \\
&\qquad\qquad\qquad\qquad - \frac{D_n^{(0)}(t',\theta)}{C_n}\{C_n - 1\}\bigg]\,dP_{\varepsilon_0,\Delta,Z}(t,1,z)
\end{aligned}
$$

$$
\begin{aligned}
(6.13) \quad &= n^{1/2} \int \{z - \hat{\eta}_n(t, \theta_0)\} \\
&\quad \times [\{D_n^{(0)}(t', \theta) - C_n d^{(0)}(t', \theta)\}] \, dP_{\varepsilon_0, \Delta, Z}(t, 1, z) + o_{p^*}(1) \\
&= \int \mathbb{G}_n[\{z - \hat{\eta}_n(t, \theta_0)\} \\
&\quad \times W\{1(\varepsilon_\theta \geq t') - d^{(0)}(t', \theta)\}] \, dP_{\varepsilon_0, \Delta, Z}(t, 1, z) + o_{p^*}(1).
\end{aligned}
$$

Similarly, we have

$$
\begin{aligned}
n^{1/2} &P[\Omega\{Z - \hat{\eta}_n(\varepsilon_0, \theta_0)\}\Delta\{\hat{\rho}_n(\varepsilon_0, \theta_0) - \rho_0(\varepsilon_0, \theta_0)\}] \\
&= \int \mathbb{G}_n[\{z - \hat{\eta}_n(t, \theta_0)\}W\{1(\varepsilon_0 \geq t) - d^{(0)}(t, \theta_0)\}] \, dP_{\varepsilon_0, \Delta, Z}(t, 1, z) \\
&\quad + o_{p^*}(1).
\end{aligned}
$$

Then, term (6.12) becomes

$$
\begin{aligned}
n^{1/2} &P[\Omega\{Z - \hat{\eta}_n(\varepsilon_0, \theta_0)\}\Delta\{\rho_0(\varepsilon_\theta, \theta) - \rho_0(\varepsilon_0, \theta_0)\}] \\
&\quad + \int \mathbb{G}_n\{z - \hat{\eta}_n(t, \theta_0)\} \\
&\qquad \times W[\{1(\varepsilon_\theta \geq t') - d^{(0)}(t', \theta)\} \\
&\qquad - \{1(\varepsilon_0 \geq t) - d^{(0)}(t, \theta_0)\}] \, dP_{\varepsilon_0, \Delta, Z}(t, 1, z) + o_{p^*}(1).
\end{aligned}
$$

Similar to the arguments following (6.10), we know that the first term above is asymptotically equivalent to $n^{1/2}(\theta - \theta_0)P[\{Z - \eta_0(\varepsilon_0, \theta_0)\}\Delta\dot{\rho}_{0\theta}(\varepsilon_0, \theta_0)]$, and the second term, above, is $o_{p^*}(1)$. So, term (6.6) can be replaced by $n^{1/2}(\theta - \theta_0)P[\{Z - \eta_0(\varepsilon_0, \theta_0)\}\Delta\dot{\rho}_{0\theta}(\varepsilon_0, \theta_0)] + o_{p^*}(1)$.

Then, from the above calculation for terms (6.5) and (6.6), we obtain

$$
\begin{aligned}
n^{1/2}&\{\Psi_n(\theta, \hat{\eta}_n(\cdot, \theta), \hat{\rho}_n(\cdot, \theta)) - \Psi_n(\theta_0, \hat{\eta}_n(\cdot, \theta_0), \hat{\rho}_n(\cdot, \theta_0))\} \\
&= -n^{1/2}(\theta - \theta_0)P\{\rho_0(\varepsilon_0, \theta_0)\dot{\eta}_{0\theta}(\varepsilon_0, \theta_0)\Delta\} \\
(6.14) \quad &\quad + n^{1/2}(\theta - \theta_0)P[\{Z - \eta_0(\varepsilon_0, \theta_0)\}\Delta\dot{\rho}_{0\theta}(\varepsilon_0, \theta_0)] + o_{p^*}(1) \\
&= n^{1/2}(\theta - \theta_0)\dot{\Psi}_\theta(\theta_0, \eta_0(\cdot, \theta_0), \rho_0(\cdot, \theta_0)) + o_{p^*}(1),
\end{aligned}
$$

which yields the asymptotic linearity (3.2) when $\theta$ is replaced by $\hat{\theta}_n$. In fact, in the above expression, we have $P[\{Z - \eta_0(\varepsilon_0, \theta_0)\}\Delta\dot{\rho}_{0\theta}(\varepsilon_0, \theta_0)] = 0$, given the equality $\eta_0(\varepsilon_0, \theta_0) = E(Z|\varepsilon_0, \Delta = 1)$, which can be verified directly (see, also, [21]). We keep it in the above calculation so as to clearly show the relationship of $\dot{\Psi}_\theta$ and $(\dot{\eta}_{0\theta}, \dot{\rho}_{0\theta})$.

Since $\hat{\theta}_n$ satisfies $\Psi_n(\hat{\theta}_n, \hat{\eta}_n(\cdot, \hat{\theta}_n), \hat{\rho}_n(\cdot, \hat{\theta}_n)) = o_{p^*}(n^{-1/2})$, showing asymptotic normality for $n^{1/2}(\hat{\theta}_n - \theta_0)$ is equivalent to showing asymptotic normality for $n^{1/2}\Psi_n(\theta_0, \hat{\eta}_n(\cdot, \theta_0), \hat{\rho}_n(\cdot, \theta_0))$. The following shows the calculation. By adding, subtracting and rearranging terms, we have

$$n^{1/2}\Psi_n(\theta_0, \hat{\eta}_n(\cdot, \theta_0), \hat{\rho}_n(\cdot, \theta_0))$$
$$= \mathbb{G}_n[\Omega\rho_0(\varepsilon_0, \theta_0)\{Z - \eta_0(\varepsilon_0, \theta_0)\}\Delta]$$

(6.15) $$\quad - \mathbb{G}_n[\Omega\hat{\rho}_n(\varepsilon_0, \theta_0)\{\hat{\eta}_n(\varepsilon_0, \theta_0) - \eta_0(\varepsilon_0, \theta_0)\}\Delta]$$

(6.16) $$\quad + \mathbb{G}_n[\Omega\{Z - \eta_0(\varepsilon_0, \theta_0)\}\Delta\{\hat{\rho}_n(\varepsilon_0, \theta_0) - \rho_0(\varepsilon_0, \theta_0)\}]$$

(6.17) $$\quad - n^{1/2}P[\Omega\hat{\rho}_n(\varepsilon_0, \theta_0)\{\hat{\eta}_n(\varepsilon_0, \theta_0) - \eta_0(\varepsilon_0, \theta_0)\}\Delta]$$

(6.18) $$\quad + n^{1/2}P[\Omega\{Z - \eta_0(\varepsilon_0, \theta_0)\}\Delta\{\hat{\rho}_n(\varepsilon_0, \theta_0) - \rho_0(\varepsilon_0, \theta_0)\}].$$

Repeatedly using similar arguments, we can show that terms (6.15) and (6.16) are $o_{p^*}(1)$. Term (6.17) can be calculated similarly, as in (6.9), but with $t = t'$, so that the lower case variable $z$ is not involved in the integrand, and $\hat{\rho}_n$ can be further replaced by $\rho_0$. Term (6.18) can be calculated similarly, as in (6.13). We then have

$$n^{1/2}\Psi_n(\theta_0, \hat{\eta}_n(\cdot, \theta_0), \hat{\rho}_n(\cdot, \theta_0))$$
$$= \mathbb{G}_n\bigg[\Omega\rho_0(\varepsilon_0, \theta_0)\{Z - \eta_0(\varepsilon_0, \theta_0)\}\Delta$$
$$\quad - \int \rho_0(t, \theta_0)d^{(0)}(t, \theta_0)^{-1}W1(\varepsilon_0 \geq t)\{Z - \eta_0(t, \theta_0)\}\, dP_{\varepsilon_0,\Delta}(t, 1)$$

(6.19) $$\quad + \int \{z - \eta_0(t, \theta_0)\}W\{1(\varepsilon_0 \geq t) - d^{(0)}(t, \theta_0)\}\, dP_{\varepsilon_0,\Delta,Z}(t, 1, z)\bigg]$$
$$\quad + o_{p^*}(1)$$
$$= \mathbb{G}_n\bigg[\Omega\rho_0(\varepsilon_0, \theta_0)\{Z - \eta_0(\varepsilon_0, \theta_0)\}\Delta$$

(6.20) $$\quad - \int \rho_0(t, \theta_0)W1(\varepsilon_0 \geq t)\{Z - \eta_0(t, \theta_0)\}\, d\Lambda_0(t)\bigg]$$
$$\quad + o_{p^*}(1),$$

which converges in distribution to a normal random variable by the central limit theorem, because the influence function in the above expression is bounded. Here, $\Lambda_0$ is the cumulative hazard function of $e_0 = T - \theta_0'Z$. So, from equation (3.2), we know that $n^{1/2}(\hat{\theta}_n - \theta_0)$ is asymptotically normal with asymptotic representation (3.3) if $\dot{\Psi}_\theta(\theta_0, \eta_0(\cdot, \theta_0), \rho_0(\cdot, \theta_0))$ is nonsingular. That the term (6.19), yielded by estimating the weight function $\rho_0(t, \theta)$, is equal to zero can be verified directly, again, by using the equality $\eta_0(\varepsilon_0, \theta_0) = E(Z|\varepsilon_0, \Delta = 1)$. Term (6.20) is obtained from the following calculation:

$$d^{(0)}(t, \theta_0) = P\{W1(Y - \theta_0'Z \geq t)\}$$

$$= P\{1(Y - \theta_0' Z \geq t)\}$$
$$= E[E\{1(Y - \theta_0' Z \geq t)|Z\}]$$
$$= E[\Pr(T - \theta_0' Z \geq t|Z)\Pr(C - \theta_0' Z \geq t|Z)]$$
$$= \int \exp\{-\Lambda_0(t)\}\{1 - G(t|z)\} \, dH(z),$$

where $G(\cdot|z)$ is the conditional distribution function of the centered censoring time $C - \theta_0' Z$ given $Z = z$, and $H$ is the marginal distribution function of covariate $Z$. On the other hand, from the joint distribution of $(\varepsilon_0, \Delta, Z)$, we obtain

$$dP_{\varepsilon_0, \Delta}(t, 1) = \left[\int \exp\{-\Lambda_0(t)\}\{1 - G(t|z)\} \, dH(z)\right] d\Lambda_0(t)$$
$$= d^{(0)}(t, \theta_0) \, d\Lambda_0(t).$$

That term (6.19) is zero becomes even more straightforward from term (6.18) if the weight function $\rho_0$ is given and, thus, need not be estimated (e.g., $\hat{\rho}_n = \rho_0 = 1$).

6.4. *Proof of Theorem 3.4.* We will sequentially show consistency, root-$n$ rate convergence and the asymptotic normality of $\hat{\theta}_n^*$. It is easy to see that $\{W(x;\alpha) : \alpha \in \mathcal{A}_0\}$ is Lipschitz in $\alpha$ and, hence, Donsker (see Example 3.2.12 of [27]), so we have that $\{\hat{\eta}_n^*\}$ and $\{\hat{\rho}_n^*\}$ are Donsker (see Section 2.10 of [27]). Based on the smoothness of $W(X;\alpha)$ in $\alpha$ and the structures of $\hat{\eta}_n$, $\hat{\rho}_n$, $\hat{\eta}_n^*$ and $\hat{\rho}_n^*$ given in (2.2), (2.3), (3.5) and (3.6), we have

$$\|W(X; \hat{\alpha}_n) - W(X; \alpha)\| \to 0, \qquad \|\hat{\eta}_n^* - \hat{\eta}_n\| \to 0 \quad \text{and} \quad \|\hat{\rho}_n^* - \hat{\rho}_n\| \to 0$$

in outer probability by the mean value theorem and boundedness of the corresponding derivatives, with respect to $\alpha$. The above three quantities are actually $O_{p^*}(n^{-1/2})$ by the root-$n$ consistency of $\hat{\alpha}_n$ and the smoothness assumption of $W(X;\alpha)$. Thus, with $\Omega_i$ replaced by $W_i$ in $\Psi_n$, we have

$$\|\Psi_n^*(\theta, \hat{\eta}_n^*(\cdot, \theta), \hat{\rho}_n^*(\cdot, \theta)) - \Psi_n(\theta, \hat{\eta}_n^*(\cdot, \theta), \hat{\rho}_n^*(\cdot, \theta))\|$$
$$(6.21) \qquad \leq \|W(X; \hat{\alpha}_n) - W(X; \alpha)\| \|\hat{\rho}_n^*(\varepsilon_\theta, \theta)\{Z - \hat{\eta}_n^*(\varepsilon_\theta, \theta)\}\Delta\|$$
$$= o_{p^*}(1)$$

by the boundedness of $\hat{\rho}_n^*(\varepsilon_\theta, \theta)\{Z - \hat{\eta}_n^*(\varepsilon_\theta, \theta)\}\Delta$. By (6.1), (6.2) and the triangle inequality, we have

$$\|\hat{\eta}_n^* - \eta_0\| \to 0 \quad \text{and} \quad \|\hat{\rho}_n^* - \rho_0\| \to 0$$

in outer probability, which by Theorem 3.1 imply that

$$\|\Psi_n(\theta, \hat{\eta}_n^*(\cdot, \theta), \hat{\rho}_n^*(\cdot, \theta)) - \Psi_n(\theta, \eta_0(\cdot, \theta), \rho_0(\cdot, \theta))\| = o_{p^*}(1),$$

since Donsker implies Glivenko–Cantelli. Hence, by the triangle inequality we have

$$\|\Psi_n^*(\theta, \hat{\eta}_n^*(\cdot, \theta), \hat{\rho}_n^*(\cdot, \theta)) - \Psi(\theta, \eta_0(\cdot, \theta), \rho_0(\cdot, \theta))\| = o_{p^*}(1),$$

which yields the consistency of $\hat{\theta}_n^*$ by the same argument as in the proof of Theorem 3.1.

From (6.21), we know that

$$\|n^{1/2}\{\Psi_n^*(\theta, \hat{\eta}_n^*(\cdot, \theta), \hat{\rho}_n^*(\cdot, \theta)) - \Psi_n(\theta, \hat{\eta}_n^*(\cdot, \theta), \hat{\rho}_n^*(\cdot, \theta))\}\|_0 = O_{p^*}(1).$$

Replacing $(\hat{\eta}_n, \hat{\rho}_n)$ with $(\hat{\eta}_n^*, \hat{\rho}_n^*)$ in (6.3), we obtain

$$\|n^{1/2}\{\Psi_n(\theta, \hat{\eta}_n^*(\cdot, \theta), \hat{\rho}_n^*(\cdot, \theta)) - \Psi(\theta, \eta_0(\cdot, \theta), \rho_0(\cdot, \theta))\}\|_0 = O_{p^*}(1).$$

Hence, by applying the triangle inequality, we have

$$\|n^{1/2}\{\Psi_n^*(\theta, \hat{\eta}_n^*(\cdot, \theta), \hat{\rho}_n^*(\cdot, \theta)) - \Psi(\theta, \eta_0(\cdot, \theta), \rho_0(\cdot, \theta))\}\|_0 = O_{p^*}(1),$$

and the same calculation as in (6.4), with $\Psi_n$ replaced by $\Psi_n^*$ and $\hat{\theta}_n$ replaced by $\hat{\theta}_n^*$, shows that $n^{1/2}(\hat{\theta}_n^* - \theta_0) = O_{p^*}(1)$.

We now prove the asymptotic normality of $n^{1/2}(\hat{\theta}_n^* - \theta_0)$. Consider the following decomposition:

$$n^{1/2}\Psi_n^*(\hat{\theta}_n^*, \hat{\eta}_n^*(\cdot, \hat{\theta}_n^*), \hat{\rho}_n^*(\cdot, \hat{\theta}_n^*))$$

(6.22) $\quad = n^{1/2}\Psi_n^*(\hat{\theta}_n^*, \hat{\eta}_n^*(\cdot, \hat{\theta}_n^*), \hat{\rho}_n^*(\cdot, \hat{\theta}_n^*)) - n^{1/2}\Psi_n(\hat{\theta}_n^*, \hat{\eta}_n(\cdot, \hat{\theta}_n^*), \hat{\rho}_n(\cdot, \hat{\theta}_n^*))$

(6.23) $\quad\quad + n^{1/2}\Psi_n(\hat{\theta}_n^*, \hat{\eta}_n(\cdot, \hat{\theta}_n^*), \hat{\rho}_n(\cdot, \hat{\theta}_n^*)) - n^{1/2}\Psi_n(\theta_0, \hat{\eta}_n(\cdot, \theta_0), \hat{\rho}_n(\cdot, \theta_0))$

(6.24) $\quad\quad + n^{1/2}\Psi_n(\theta_0, \hat{\eta}_n(\cdot, \theta_0), \hat{\rho}_n(\cdot, \theta_0)) - n^{1/2}\Psi_n(\hat{\theta}_n, \hat{\eta}_n(\cdot, \hat{\theta}_n), \hat{\rho}_n(\cdot, \hat{\theta}_n))$

(6.25) $\quad\quad + n^{1/2}\Psi_n(\hat{\theta}_n, \hat{\eta}_n(\cdot, \hat{\theta}_n), \hat{\rho}_n(\cdot, \hat{\theta}_n)).$

Then, applying (6.14) to (6.23) and (6.24), respectively, we can replace (6.23) with

(6.26) $\quad\quad n^{1/2}(\hat{\theta}_n^* - \theta_0)\dot{\Psi}_\theta(\theta_0, \eta_0(\cdot, \theta_0), \rho_0(\cdot, \theta_0)) + o_{p^*}(1)$

and replace (6.24) with

(6.27) $\quad\quad -n^{1/2}(\hat{\theta}_n - \theta_0)\dot{\Psi}_\theta(\theta_0, \eta_0(\cdot, \theta_0), \rho_0(\cdot, \theta_0)) + o_{p^*}(1).$

Term (6.25), clearly, is $o_{p^*}(1)$. We then calculate term (6.22). Let

$$\hat{\eta}_{n,\alpha}(t, \theta) = \mathbb{P}_n\{W(X; \alpha)1(\varepsilon_\theta \geq t)Z\}/\mathbb{P}_n\{W(X; \alpha)1(\varepsilon_\theta \geq t)\},$$

$$\hat{\rho}_{n,\alpha}(t, \theta) = \mathbb{P}_n\{W(X; \alpha)1(\varepsilon_\theta \geq t)\}/\mathbb{P}_n\{W(X; \alpha)\}.$$

Then, we have $\hat{\eta}_n \equiv \hat{\eta}_{n,\alpha_0}$, $\hat{\rho}_n \equiv \hat{\rho}_{n,\alpha_0}$, $\hat{\eta}_n^* \equiv \hat{\eta}_{n,\hat{\alpha}_n}$, and $\hat{\rho}_n^* \equiv \hat{\rho}_{n,\hat{\alpha}_n}$. Let

$$\Phi_n(\alpha, \theta) = \mathbb{P}_n[W(X; \alpha)\hat{\rho}_{n,\alpha}(\varepsilon_\theta, \theta)\{Z - \hat{\eta}_{n,\alpha}(\varepsilon_\theta, \theta)\}\Delta].$$

It can be seen by direct calculation that the second derivative of $\Phi_n(\alpha, \theta)$ to $\alpha$ is bounded with outer probability 1. So, by the Taylor expansion, we have

$$\Psi_n^*(\theta, \hat{\eta}_n^*, \hat{\rho}_n^*) - \Psi_n(\theta, \hat{\eta}_n, \hat{\rho}_n) = \Phi_n(\hat{\alpha}_n, \theta) - \Phi_n(\alpha_0, \theta)$$
$$= \dot{\Phi}_{n,\alpha}(\alpha_0, \theta)(\hat{\alpha}_n - \alpha_0) + o_{p^*}(n^{-1/2}),$$

where

$$\dot{\Phi}_{n,\alpha}(\alpha_0, \theta) = \mathbb{P}_n \bigg[ \hat{\rho}_{n,\alpha_0}(\varepsilon_\theta, \theta) \{ Z - \hat{\eta}_{n,\alpha_0}(\varepsilon_\theta, \theta) \} \frac{\partial W(X; \alpha)}{\partial \alpha'} \bigg|_{\alpha=\alpha_0} \Delta$$
$$+ W(X; \alpha_0) \{ Z - \hat{\eta}_{n,\alpha_0}(\varepsilon_\theta, \theta) \} \frac{\partial \hat{\rho}_{n,\alpha}(\varepsilon_\theta, \theta)}{\partial \alpha'} \bigg|_{\alpha=\alpha_0} \Delta$$
$$+ W(X; \alpha_0) \hat{\rho}_{n,\alpha_0}(\varepsilon_\theta, \theta) \bigg\{ -\frac{\partial \hat{\eta}_{n,\alpha}(\varepsilon_\theta, \theta)}{\partial \alpha'} \bigg\}_{\alpha=\alpha_0} \Delta \bigg].$$

It is also easy to see, by direct calculation, that $\{\partial \hat{\eta}_{n,\alpha}/\partial \alpha|_{\alpha=\alpha_0} : \theta \in \Theta_0\}$ and $\{\partial \hat{\rho}_{n,\alpha}/\partial \alpha|_{\alpha=\alpha_0} : \theta \in \Theta_0\}$ are (componentwise) Glivenko–Cantelli, so, with outer probability 1, we have

$$\dot{\Phi}_{n,\alpha}(\alpha_0, \hat{\theta}_n^*) \to P[\rho_0(\varepsilon_0, \theta_0) \{ Z - \eta_0(\varepsilon_0, \theta_0) \} (\dot{W}_\alpha(X; \alpha_0))' \Delta]$$
$$+ P[W(X; \alpha_0) \{ Z - \eta_0(\varepsilon_0, \theta_0) \} A_1(\varepsilon_0, \theta_0) \Delta]$$
(6.28)
$$- P[W(X; \alpha_0) \rho_0(\varepsilon_0, \theta_0) A_2(\varepsilon_0, \theta_0) \Delta]$$
$$= P[\rho_0(\varepsilon_0, \theta_0) \{ Z - \eta_0(\varepsilon_0, \theta_0) \} (\dot{W}_\alpha(X; \alpha_0))' \Delta]$$
$$- P[\rho_0(\varepsilon_0, \theta_0) A_2(\varepsilon_0, \theta_0) \Delta]$$
$$\equiv -B,$$

where $\dot{W}_\alpha(X; \alpha) = \partial W(X; \alpha)/\partial \alpha$, $A_1$ is the limit of $\partial \hat{\rho}_{n,\alpha}/\partial \alpha'|_{\alpha=\alpha_0, \theta=\hat{\theta}_n^*}$ and $A_2$ is the limit of $\partial \hat{\eta}_{n,\alpha}/\partial \alpha'|_{\alpha=\alpha_0, \theta=\hat{\theta}_n^*}$. The term (6.28) is zero since $E(Z|\varepsilon_0, \Delta = 1) = \eta_0(\varepsilon_0, \theta_0)$. Note that $E(W|X) = 1$ is also used in the above calculation. It can be directly verified that

$$A_2(t, \theta_0) = \frac{1}{P\{(1(\varepsilon_0 \geq t)\}}[P\{1(\varepsilon_0 \geq t) Z (\dot{W}_\alpha(X; \alpha_0))'\}$$
$$- \eta_0(t, \theta_0) P\{(\dot{W}_\alpha(X; \alpha_0))' 1(\varepsilon_0 \geq t)\}].$$

Hence, we have

(6.29) $\quad \Psi_n^*(\theta, \hat{\eta}_n^*, \hat{\rho}_n^*) - \Psi_n(\theta, \hat{\eta}_n, \hat{\rho}_n) = -B(\hat{\alpha}_n - \alpha_0) + o_{p^*}(n^{-1/2}).$

Replacing (6.22), (6.23) and (6.24) by (6.29), (6.26) and (6.27), respectively, we obtain

$$n^{1/2}(\hat{\theta}_n^* - \theta_0) = n^{1/2}(\hat{\theta}_n - \theta_0) + \{\dot{\Psi}_\theta(\theta_0, \eta_0, \rho_0)\}^{-1} B n^{1/2}(\hat{\alpha}_n - \alpha_0) + o_{p^*}(1).$$

26  B. NAN, J. D. KALBFLEISCH AND M. YU

By (3.3) we know that $\hat{\theta}_n$ is an asymptotically linear estimator. Given that $\hat{\alpha}_n$ is also an asymptotically linear estimator, we know that $n^{1/2}(\hat{\theta}_n - \theta_0)$ and $n^{1/2}(\hat{\alpha}_n - \alpha_0)$ are asymptotically jointly normal by the multivariate central limit theorem. Hence, by [18], we know that $n^{1/2}(\hat{\theta}_n^* - \theta_0)$ is asymptotically normal with variance given in (3.7).

B. Nan
J. D. Kalbfleisch
Department of Biostatistics
University of Michigan
1420 Washington Heights
Ann Arbor, Michigan 48109-2029
USA
E-mail: bnan@umich.edu
jdkalbfl@umich.edu

M. Yu
Department of Medicine and Biostatistics
Indiana University School of Medicine
535 Barnhill Drive, RT 380F
Indianapolis, Indiana 46202
USA
E-mail: meyu@iupui.edu